\documentclass[final]{ifacconf}
\PassOptionsToPackage{dvipsnames}{xcolor}

\usepackage{caption}

\usepackage[ieee]{mathematix}
\usepackage{units}
\usepackage{xspace}

\newcommand{\OpEn}{\texttt{OpEn}\xspace}
\newcommand{\C}{\texttt{C}\xspace}
\newcommand{\Cpp}{\texttt{C++}\xspace}
\newcommand{\opengen}{\texttt{opengen}\xspace}

\newcommand{\ubarnu}{{u}^{\nu}}
\newcommand{\ubarnuplus}{{u}^{\nu{+}1}}
\newcommand{\ubarhalf}{{u}^{\nu{+}\nicefrac{1}{2}}}

\usepackage{xcolor}

\usepackage{listings}
\usepackage{textcomp}
\lstdefinestyle{myMatlab}{
    language=Matlab,                        
    upquote=true,
    frame=single,                           
    basicstyle=\small\ttfamily,             
    backgroundcolor=\color{CornflowerBlue!8},
    keywordstyle=[1]\color{NavyBlue}\bfseries,  
    keywordstyle=[2]\color{purple},         
    keywordstyle=[3]\color{black}\bfseries, 
    identifierstyle=,                       
    commentstyle=\usefont{T1}{pcr}{m}{sl}\color{CommentGreen}\small,
    stringstyle=\color{purple},             
    showstringspaces=false,                 
    tabsize=5,                              
    morekeywords={properties,methods,classdef},
    morekeywords=[2]{handle,syms},
    morecomment=[l][\color{blue}]{...},     
    numbers=none,                           
    firstnumber=1,                          
    numberstyle=\tiny\color{blue},          
    stepnumber=1,                            
    xleftmargin=10pt,
    xrightmargin=10pt
}

\lstdefinestyle{myPython}{
    language=Python,                        
    upquote=true,
    basicstyle=\footnotesize\ttfamily,             
    backgroundcolor=\color{gray!10},
    keywordstyle=[1]\color{NavyBlue},  
    keywordstyle=[2]\color{ForestGreen},         
    keywordstyle=[3]\color{magenta}, 
    identifierstyle=,                       
    commentstyle=\usefont{T1}{pcr}{m}{sl}\color{MidnightBlue}\small,
    stringstyle=\color{purple},             
    showstringspaces=false,                 
    tabsize=5,                              
    morekeywords={import,as},
    morekeywords=[2]{og,ogc,cs},
    morekeywords=[3]{radius},
    morecomment=[l][\color{blue}]{...},     
    numbers=none,                           
    firstnumber=1,                          
    numberstyle=\tiny\color{blue},          
    stepnumber=1,                            
    xleftmargin=3pt,
    xrightmargin=3pt
}

  \usetikzlibrary{matrix}
  \usetikzlibrary{shapes,arrows}
  \usetikzlibrary{calc}
  \usetikzlibrary{patterns}
  \usetikzlibrary{intersections}
  \usetikzlibrary{decorations.markings}

\begin{document}

\begin{frontmatter}

\title{OpEn: Code Generation for Embedded Nonconvex Optimization%
\thanksref{footnoteinfo}} 

%
%
\thanks[footnoteinfo]{The work of P. Patrinos was supported by: FWO projects: No. G086318N; No. G086518N; 
        Fonds de la Recherche Scientifique -- FNRS, the Fonds Wetenschappelijk Onderzoek 
        -- Vlaanderen under EOS Project No. 30468160 (SeLMA), 
        Research Council KU Leuven C1 project No. C14/18/068 and the 
        Ford--KU Leuven Research Alliance project. Corresponding author email: \texttt{p.sopasakis@qub.ac.uk}.}

\author[QUB]{Pantelis Sopasakis},
\author[WIDEFIND]{Emil Fresk}
and
\author[KUL]{Panagiotis Patrinos}

\address[QUB]{Queen's University Belfast, School of Electronics, Electrical Engineering and Computer Science
(EEECS) and Centre for Intelligent and Autonomous Manufacturing Systems (i-AMS), Ashby Building, Stranmillis Road,
BT9 5AG Belfast, Northern Ireland}
\address[WIDEFIND]{WideFind AB, Aurorum 1C, 977 75 Lule\aa, Sweden}
\address[KUL]{KU Leuven, Department of Electrical Engineering (ESAT),
              STADIUS Center for Dynamical Systems, Signal Processing and Data Analytics
              \& Optimization in Engineering (OPTEC),
              Kasteelpark Arenberg 10, 3001 Leuven, Belgium.}


\begin{abstract}  
We present Optimization Engine (\OpEn): an open-source code generation tool for real-time embedded
nonconvex optimization, which implements  
a novel numerical method. \OpEn combines the proximal averaged Newton-type method for optimal 
control (PANOC) with the penalty and augmented Lagrangian methods to compute approximate 
stationary points of nonconvex problems. The proposed method involves very simple algebraic
operations such as vector products, has a low memory footprint and exhibits very good convergence
properties that allow the solution of nonconvex problems on embedded devices. 
\OpEn's core solver is written is Rust --- a modern, high-performance, memory-safe and thread-safe systems 
programming language --- while users can call it from Python, MATLAB, \C, \Cpp or over a TCP socket.
\end{abstract}

\end{frontmatter}

\section{Introduction}

In embedded applications it is often necessary to solve optimization problems
in real time. Typical examples involve moving horizon estimation, 
model predictive control and online learning \citep{Ferreau+2017}. 
The distinctive features of embedded real time
computing are the limited computing and storage capabilities of hardware devices,
the presence of stringent runtime requirements,
the increased need for memory safety and reliability,
the need for simple numerical methods which can be verified easily,
and the need for design paradigms that involve code generation.

A significant research effort has been dedicated to numerical optimization software
with code generation capabilities for convex problems such as 
CVXGen \citep{Mattingley2012,MattingleyBoyd2010b},
SPLIT \citep{SHUKLA201714386}, 
OSQP \citep{osqp-codegen}, and 
$\mu$AO-MPC \citep{muAOMPC2013}.
\cite{KouzoupisBenchmarks2015} provide benchmarks of numerical optimization 
methods for embedded quadratic programming.

Nonconvex optimization problems are typically solved with 
sequential convex programming (SCP) \citep{Zillober2004} and interior point (IP) methods 
\citep{Zanelli+2017,GopalBiegler1998}. Several software are available for 
nonconvex optimization such as SNOPT \citep{Gill2005}, Acado \citep{Verschueren+2018,Houska2011a,Houska2011b}
and IPOPT \citep{Wachter2006}.

Nonetheless, SCP and IP methods have a high per-iteration computation cost and involve expensive 
operations such as the solution of quadratic programming problems or the solution of 
linear systems. Was it not for their slow convergence, first-order methods such as 
the projected gradient method, would be ideal for embedded applications. 
\cite{panoc2017} proposed the proximal averaged Newton-type method for optimal 
control (PANOC), which uses exactly the same oracle as the projected gradient method
but exhibits superior convergence characteristics and has been found to outperform
SCP and IP methods \citep{agv2018,mav2019}.

PANOC can solve problems with a smooth cost function and simple constraints on 
which one can compute projections. More general constraints can be relaxed and 
replaced by penalty functions (soft constraints). The penalty method can then 
be applied to obtain near-feasible solutions \citep{Hermans:IFAC:2018}.
However, as the penalty parameter needs to grow unbounded, it is likely that it 
will overshadow the cost function and lead to highly ill conditioned problems.
In this paper we propose a numerical method to determine stationary points of 
constrained nonconvex optimization problems by combining the penalty and the 
augmented Lagrangian methods.
Moreover, we present \OpEn: an open-source software an open-source code generation tool for
real-time embedded nonconvex optimization with guaranteed safe memory management.
In Section \ref{sec:simulations} we present two simulation examples -- a model predictive 
controller for obstacle avoidance and a constrained nonlinear estimator --
and provide comparisons with IPOPT and \texttt{scipy}'s SQP solver.

\subsection*{Preliminaries and Notation}
We define $\N_{N}=\{0,1,\ldots, N\}$.
We denote the set of extended-real numbers by $\barre=\R\cup\{\infty\}$.
The Jacobian of a mapping $F:\R^m\to\R^n$ is denoted by $JF$.
For two vectors, $x, x'\in\R^n$, the relation $x\leq x'$ is meant 
in the pointwise sense.
For $x\in\R$, we define $\plus{x}=\max\{0, x\}$ and for $x\in\R^n$, $\plus{x}$
is defined element wise.
The projection of a vector $x\in\R^n$ on a nonempty, closed, convex set $C\subseteq\R^n$ is defined as 
\(\proj_{C}(x) = \argmin_{y\in C}\|y-x\|\) and the distance of $x\in\R^n$ from $C$ is defined
as $\dist_C(x) = \min_{y\in C}\|y-x\|$.
The support function of a nonempty, closed, convex set $C\subseteq\R^n$ is a function $\delta_C^*:\R^n\to\barre$
defined as $\delta_C^*(y) = \sup_{x\in C}y^\top x$. 
For a cone $K\subseteq\R^n$, we define its polar to be the cone $K^\circ = \{y\in\R^n{}:{}y^\top x \leq 0\}$.
The normal cone of a set $C\subseteq\R^n$ is defined as that subdifferential of $\delta_C$, that is,
$N_C(x)=\partial \delta_C(x)$, for $x\in\R^n$.

\section{Problem statement}
Consider the following parametric optimization problem
\begin{subequations}\label{eq:original_problem}
\begin{align}
 \mathbb{P}(p){}:{}\minimize_{u{}\in{}U}\,&f(u, p),
 \\
 \subjto\,& F_1(u, p) \in C,\label{eq:constraints:alm}
 \\
 &F_2(u, p) = 0,\label{eq:constraints:pm}
\end{align}
\end{subequations}
where for each $p\in\R^{n_p}$, $f({}\cdot{}, p): \R^{n} \to \R$ is a (possibly nonconvex)
continuously differentiable function with $L_f$-Lipschitz gradient and
$U\subseteq\R^{n}$ is a nonempty, closed --- but not necessarily convex ---
set such that we can compute projections on it.
Moreover, $F_1({}\cdot{}, p):\R^{n}\to\R^{n_1}$ is a smooth mapping with Lipschitz-continuous
Jacobian which is bounded on $U$.
Set $C\subseteq\R^{n_1}$ is a convex set from which we can compute distances.
Mapping $F_2({\cdot}, p):\R^{n}\to\R^{n_2}$ is such that
$\|F_2({}\cdot{}; p)\|^2$ is a continuously differentiable function with
Lipschitz-continuous gradient.

In particular, constraints \eqref{eq:constraints:alm} will be treated using the proposed
augmented Lagrangian method, while constraints \eqref{eq:constraints:pm} will be accounted
for using the penalty method as we shall discuss in Section \ref{sec:numerical}.
Constraints of both types can be present at the same time.

Constraints \eqref{eq:constraints:pm} can accommodate multiple vertical complementarity
constraints \citep{scheel2000complem,Izmailov2012} of the form
\begin{equation}
 \prod_{i=1}^{m}\plus{h_i(u, p)} {}={} 0,
\end{equation}
a special case of which arises in collision avoidance problems
\citep{Hermans:IFAC:2018,agv2018,mav2019}. Such constraints do not satisfy the smoothness
conditions, so they cannot be described in terms of the ALM-type constraints in
\eqref{eq:constraints:alm}.

Set $U$ can be a Euclidean ball of radius $r$ centered at a point $u_c$,
$\mathcal{B}_2(u_c, r)=\{u\in\R^{n}{}:{}\|u-u_c\|_2\leq r\}$, an infinity ball,
$\mathcal{B}_\infty(u_c, r)=\{u\in\R^{n}{}:{}\|u-u_c\|_\infty\leq r\}$,
a rectangle, $\{u\in\R^n{}:{} u_{\mathrm{min}} \leq u \leq u_{\mathrm{max}}\}$ (where some coordinates of
$u_{\mathrm{min}}$ and $u_{\mathrm{max}}$ are allowed to be equal to $-\infty$ and $\infty$
respectively), a finite set
$\{u^{(1)},u^{(2)},\ldots,u^{(q)}\}$, or a Cartesian product of any such
sets.

The above formulation can accommodate equality constraints, either through constraint
\eqref{eq:constraints:pm}, or by using $C=\{0\}$. It can also accommodate inequality constraints
of the general form
\begin{equation}
 H(u, p) \leq 0,
\end{equation}
by using either $F_1(u, p) = H(u, p)$ and
\(
	C
  {}={}
	\{
	    v\in\R^{n_1}
	  {}:{}
	    v \leq 0
	\}
\),
or
\(
	F_2(u, p)
  {}={}
	\plus{{}H(u, p){}}
\).
Furthermore, constraint \eqref{eq:constraints:alm} allow the designer to encode norm
bounds using $C=\mathcal{B}_2(u_c, r)$ or $C=\mathcal{B}_\infty(u_c, r)$, or constraints
involving second-order cones,
\(
	C
  {}={}
	\mathrm{SOC}_\alpha^{n_1}
  {}\dfn{}
	\{
	    u=(x, t),
	    x\in\R^{n_1-1},
	    t\in\R:
	    \|x\| \leq \alpha t
	\}
\), for some \(\alpha > 0\).

Optimal control and estimation problems can be written in the form of Problem
\eqref{eq:original_problem} in a number of different ways.
For instance, consider the following discrete-time optimal control problem
\begin{subequations}
 \begin{align}
  \mathbb{P}(x_0):\minimize\,&\ell_N(x_N) + \sum_{t=0}^{N-1}\ell_t(x_t, u_t)
  \\
  \subjto\,&u_t \in U_t,
  \\
  &x_{t+1} = \Phi_t(x_t, u_t),\label{eq:optimal_control:dynamics}
  \\
  &H_t(x_t, u_t) \leq 0,\label{eq:optimal_control:xu_constraints}
  \\
  &H_N(x_N) \leq 0,\label{eq:optimal_control:terminal_constraints}
 \end{align}
\end{subequations}
for \(t\in\N_{N-1}\), where $U_t$ are sets on which one can easily compute projections.
One may follow the \textit{single shooting} approach, in which the minimization is carried out
over the sequence of control actions, $u=(u_0, \ldots, u_{N-1})$, and the sequence of states is eliminated
by defining a sequence of function $F_0(u, x_0) = x_0$ and $F_{t+1}(u, x_0) = \Phi_t(F_t(u, x_0), u_t)$
for \(t\in\N_{N-1}\). We may then define the cost function
\begin{equation}\label{eq:single_shooting}
 f(u, x_0) = \ell_N(F_N(u, x_0)) + \sum_{t=0}^{N-1}\ell_t(F_t(u, x_0), u_t).
\end{equation}
A second alternative is to follow the \textit{multiple shooting} approach where the minimization
is carried out over $u=(u_0, \ldots, u_{N-1},x_1,\ldots, x_{N})$ and the system dynamics
in Equation \eqref{eq:optimal_control:dynamics} are treated as constraints.

Joint state-input constraints of the form \eqref{eq:optimal_control:xu_constraints} and
terminal constraints of the form \eqref{eq:optimal_control:terminal_constraints} can be
cast using either $F_1$ or $F_2$ as discussed above.

\section{Numerical algorithm}\label{sec:numerical}

\subsection{Augmented Lagrangian and Penalty Methods}
In order to solve the original optimization problem, we shall use the augmented Lagrangian
method \citep{BirginMartinez2014}.
Hereafter, we shall drop the parameter $p$ for the sake of simplicity.
We introduce the \textit{augmented Lagrangian function}
\begin{multline}
 L_c(u, v, y) \dfn f(u) + \<y, F_1(u)-v\> \\+ \tfrac{c}{2}\|F_1(u)-v\|^2 +  \tfrac{c}{2}\|F_2(u)\|^2,
\end{multline}
defined for $u\in U$ and $v \in C$, where $c>0$ is a constant penalty parameter
and $y\in\R^{n_1}$ is the vector of Lagrange multipliers for constraints \eqref{eq:constraints:alm}.
We can then show that

\begin{prop}
The following holds for all $y\in\R^{n_1}$ and $c>0$
\begin{multline}
 \min_{u\in U,v\in C}L_c(u, v, y) = -\tfrac{1}{2c}\|y\|^2
 + \min_{u\in U} \psi(u; c, y),
\end{multline}
where
\begin{equation}
 \psi(u; c, y) {\dfn} f(u) \,{+}\, \tfrac{c}{2}\Big[\hspace{-0.2em}\dist_C^2(F_1(u){+}\tfrac{1}{c} y)
 \,{+}\, \|F_2(u)\|^2\Big]\!\,.\hspace{-0.2em}
\end{equation}
\end{prop}
\begin{pf}
 For $u{}\in{}U$ and $v{}\in{}C$, we have that
 \begin{align*}
  {}&{}\<y, F_1(u)-v\> + \tfrac{c}{2}\|F_1(u)-v\|^2
  \\
  {}={}&c\<\tfrac{1}{c}y, F_1(u){-}v\> + \tfrac{c}{2}\|F_1(u){-}v\|^2 + \tfrac{1}{2c}\|y\|^2 -  \tfrac{1}{2c}\|y\|^2
  \\
  {}={}&\tfrac{c}{2}\|F_1(u)-v+\tfrac{1}{c}y\|^2 - \tfrac{1}{2c}\|y\|^2,
 \end{align*}
 therefore,
 \begin{multline}
 \min_{u\in U,v\in C}L_c(u, v, y) = - \tfrac{1}{2c}\|y\|^2 + \min_{u\in U} \Big\{ f(u)
 \\+ \tfrac{c}{2}\|F_2(u)\|^2
 + \smashunderbracket{\min_{v\in C}\tfrac{c}{2}\|F_1(u)-v+\tfrac{1}{c}y\|^2}{\dist_C^2(F_1(u)+c^{-1}y)}\Big\},
\end{multline}
\vspace{0.3em}

which proves the assertion.~\hfill{$\Box$}
\end{pf}

The above choice of an augmented Lagrangian function leads to Algorithm
\ref{alg:alm_pm} where the \textit{inner problem} takes the form
\begin{equation}
 \mathbb{P}_{\mathrm{in}}(c, y){}:{} \minimize_{u\in U}\psi(u; c, y).
\end{equation}
This problem has a smooth cost function and simple constraints, therefore, it
can be solved using PANOC as we discuss in Section \ref{sec:inner_problems}.
The algorithm updates both a penalty parameter, $c$, as well as a vector of Lagrange multipliers,
$y\in\R^{n_1}$, corresponding to constraint \eqref{eq:constraints:alm} according to
\begin{equation}
	y^{\nu+1}
 {}={}
	\bar{y}^{\nu} + c_\nu(F_1(u^{\nu+1})
	 - \proj_C(F_1(u^{\nu+1}) + c^{-1}_\nu\bar{y}^\nu)).
\end{equation}

\begin{algorithm}[h]
\caption{Augmented Lagrangian and penalty method}\label{alg:alm_pm}
 \begin{algorithmic}[1]
  \Require \(u^0\in\R^{n_u}\) (initial guess), 
           \(p\in\R^{n_p}\) (parameter),
           \(y^0\in\R^{n_1}\) (initial guess for the Lagrange multipliers), 
           \(\epsilon, \delta > 0\) (tolerances),
           \(\beta\) (tolerance decrease coefficient), 
           \(\rho\) (penalty update coefficient), 
           \(\theta\) (sufficient decrease coefficient),
           \(Y\subseteq \dom \delta^*_C\) (compact set)
  \Ensure \((\epsilon, \delta)\)-approximate solution $(u^\star, y^\star)$
  \State $\bar{\epsilon}_0 = \epsilon_0$
  \For{$\nu=0,\ldots, \nu_{\mathrm{max}}$}
    \State $\bar{y}^{\nu}{}={}\proj_Y(y^\nu)$
    \State $u^{\nu+1}$ is a solution of $\mathbb{P}_{\mathrm{in}}(c_\nu, \bar{y}^\nu)$ with tolerance $\bar\epsilon$ and initial guess $u^\nu$,
           using Algorithm \ref{alg:panoc}
    \State $y^{\nu+1} {}={} \bar{y}^{\nu} + c_\nu(F_1(u^{\nu+1}) - \proj_C(F_1(u^{\nu+1}) + c^{-1}_\nu\bar{y}^\nu))$
    \State $z_{\nu+1} {}={} \Vert y^{\nu+1} - \bar{y}^\nu \Vert_\infty$, 
           $t_{\nu+1} = \Vert F_2(u^{\nu+1}) \Vert_\infty$
    \If{$z_{\nu+1} \leq c_\nu \delta$, $t_{\nu+1} \leq \delta$ \textbf{and} $\bar{\epsilon}_\nu \leq \epsilon$}
      \State \textbf{return} $(u^\star, y^\star)=(u^{\nu+1}, y^{\nu+1})$
    \ElsIf{$\nu>0$, $z_{\nu+1} > \theta z_{\nu}$, $t_{\nu+1} > \theta t_\nu$}
      \State $c_{\nu+1} = \rho{}c_{\nu}$
    \EndIf{}
    \State $\bar\epsilon_{\nu+1} = \beta\bar{\epsilon}_{\nu}$
  \EndFor{}
\end{algorithmic}
\end{algorithm}

The most critical tuning parameter of the algorithm is the \textit{penalty update factor},
$\rho>1$. Typical values are between $1.1$ and $20$. Low values lead to sequences of inner
problems which are ``similar'' enough so that $u^\nu$ is a good initial guess, however,
at the expense of slow convergence of the infeasibility as quantified by $z_{\nu}$ and $t_{\nu}$.
On the other hand, larger values of $\rho$ create more dissimilar inner problems thus
deeming $u^{\nu}$ potentially not too good initial guesses, but lead to fewer outer iterations.

Note that if $u^{\nu+1}$ and $y^{\nu+1}$ lead to a significant improvement in infeasibility
as measured by $z_{\nu+1}$ and $t_{\nu+1}$, then the penalty parameter is not increased.
Moreover, the tolerance $\epsilon$ of the inner problems can be relaxed by starting from a
given initial tolerance $\epsilon_0$ and decreasing it gradually until it reaches the target
inner tolerance, $\epsilon$.

Set $Y$ is taken to be a compact subset of $\dom \delta_C^*$.
If $C$ is the Minkowski sum of a compact set and a cone $K\subseteq\R^n$,
then $\dom \delta_C^* = K^\circ$.
For example, if $C$ is compact, we may select $Y = [-M, M]^{n_1}$, whereas
if $C=\R_+^{n_1}$, we select $Y=[-M, 0]^{n_1}$.

\subsection{Inner Problems}\label{sec:inner_problems}
Problem $\mathbb{P}_{\mathrm{in}}(c, y)$ has a smooth cost function with gradient
\begin{multline}
 \nabla_u \psi(u; c, y) = \nabla f(u) + c JF_1(u)^\top \Big[F_1(u) + c^{-1}y\\
 - \proj_C(F_1(u) + c^{-1}y)\Big] + \tfrac{c}{2} \nabla \|F_2(u)\|^2.
\end{multline}
Under the assumptions on $F_1$ and $F_2$, \(\nabla_u \psi\) is Lipschitz in
$u$ with some Lipschitz modulus $L_\psi>0$, therefore, $\mathbb{P}_{\mathrm{in}}$ falls into
the framework of PANOC \citep{panoc2017}. This gradient can be computed by automatic
differentiation software such as CasADi \citep{Andersson2019}.

Our aim is to determine a $u^\star$ that satisfies the first-order optimality
conditions of $\mathbb{P}_{\mathrm{in}}(c, y)$
\begin{equation}\label{eq:simple_problem:optimality_conditions}
 u^\star = T_\gamma(u^\star; c, y),
\end{equation}
where $T_\gamma:\R^{n}{}\to{}U$, with $\gamma>0$, is the projected gradient operator,
$T_\gamma(u; c, y) = \proj_U(u - \gamma\nabla_u \psi(u; c, y))$.
To that end, we may use the proximal gradient method,
\begin{equation}
 u^{\nu+1} = T_\gamma(u^\nu; c, y),
\end{equation}
the accumulation points of which satisfy the fixed point condition \eqref{eq:simple_problem:optimality_conditions}
whenever $\gamma < \nicefrac{2}{L_\psi}$ \cite[Prop. 2.3.2]{nocedal2006numerical}. The iterations are very simple, however, convergence can be
particularly slow, especially for ill-conditioned problems.

Instead, we follow the approach of PANOC \citep{panoc2017}: we define the \textit{fixed point residual}
operator $R_\gamma:\R^{n}\to\R^{n}$ as
\begin{equation}
 R_\gamma(u; c, y) = \gamma^{-1}(u - T_\gamma(u; c, y)).
\end{equation}
Then a $u^\star\in\R^{n}$ satisfies \eqref{eq:simple_problem:optimality_conditions} if and only
if it is a zero of $R_\gamma$, that is, if it solves the nonlinear equation
\(
	  R_\gamma(u^\star; c, y)
      {}={}
	  0
\).
The main idea is that thereon we can apply a Newton-type method of the form
\(
	u^{\nu+1}
    {}={}
	u^{\nu} -H_\nu R_\gamma(u^\nu; c, y),
\)
where $H_\nu:\R^{n}\to\R^{n}$ are invertible linear operators
that satisfy the secant condition
\(
	u^{\nu+1} - u^{\nu}
    {}={}
	H_\nu (R_\gamma(u^{\nu+1}; c, y) - R_\gamma(u^{\nu}; c, y)).
\)
This motivates the use of quasi-Newtonian directions such as L-BFGS,
which is known to yield a good convergence rate at a low memory footprint
\cite[Sec.~7.2]{nocedal2006numerical}.
However, convergence can only be guaranteed when $u^0$ is in a neighborhood
of a critical point $u^\star$. To overcome this limitation, PANOC employs the
\textit{forward-backward envelope}: a real-valued, continuous merit function
$\varphi_\gamma:\R^{n}\to\R$ for Problem $\mathbb{P}_{\mathrm{in}}(c, y)$ \citep{themelis2016forward}.
For $0< \gamma < \nicefrac{1}{L_\psi}$, $\varphi_\gamma$ shares the same (local/strong)
minima with $\mathbb{P}_{\mathrm{in}}$ \citep{panoc2017}. Function $\varphi_\gamma$
is given by
\begin{multline*}
 \varphi_\gamma(u; c, y) = \psi(u; c, y) \,{-}\, \tfrac{\gamma}{2}\|\nabla \psi(u; c, y)\|^2
 \\
 {}+{} \tfrac{1}{2\gamma}\dist_U^2(u - \gamma \nabla \psi(u; c, y)).
\end{multline*}
Note that \(\varphi_\gamma(u; c, y)\) is computed at the cost of one projected gradient
step on $\psi$.

PANOC is shown in Algorithm \ref{alg:panoc}. PANOC combines \textit{safe}
projected gradient updates with \textit{fast} quasi-Newtonian directions, which are computed
with L-BFGS (see line \ref{alg:panoc::update} in Algorithm \ref{alg:panoc}) based on
the modification of L-BFGS proposed by \cite{Li2001}. With this choice of averaged directions
and under mild conditions, eventually only fast quasi-Newtonian directions will be activated
leading to very fast convergence \citep{panoc2017}.

PANOC uses the same oracle as the projected gradient method --- it only requires $\psi$, $\nabla \psi$
and $\proj_U$ --- and involves a very simple decrease criterion on $\varphi_\gamma$.
The computation of the L-BFGS directions, $d^\nu$, requires merely $4{}\mu{}n$, where $\mu$
is the L-BFGS memory.

\begin{algorithm}
\caption{PANOC method for solving the inner problem, \(\mathbb{P}_{\mathrm{in}}(c, y)\)}\label{alg:panoc}
 \begin{algorithmic}[1]
  \Require \(u^0\in\R^{n}\) (initial guess),\
 	\(L_{\psi,0}>0\) (estimate of the Lipschitz constant of \(\nabla \psi({}\cdot{}; c, y)\)),\
 	\(\mu\) (memory length of L-BFGS),
 	\(\epsilon>0\) (tolerance),
 	\(\nu_{\mathrm{max}}\) (max. iterations)%
   \Ensure Approximate solution \(u^{\star}\)
   \State \(L {}\gets{} L_{\psi, 0}\)
   \State Choose \(\gamma\in(0,\nicefrac{1}{L})\), \(\sigma\in(0,\tfrac{\gamma}{2}(1-\gamma L))\)
   \For{\(\nu=0,1,\ldots,\nu_{\mathrm{max}}\)}
 	\State
 		Compute \(\nabla \psi({u}^\nu; c, y)\) with automatic differentiation
	\State\label{alg:panoc::FB}
		\(
			\ubarhalf
		{}\gets{}
			T_{\gamma}(\ubarnu; c, y)
		\)
	\State\label{alg:panoc::r}
		\(
			r^{\nu}
		    {}\leftarrow{}
			u^\nu{}-{}\ubarhalf
		\)
	\State \algorithmicif{} 
		{\(
		      \|
			\gamma^{-1} r^{\nu} 
		        {+} 
		        \nabla \psi(\ubarhalf; c, y) 
		        {-} 
		        \nabla \psi(\ubarnu; c, y) 
		      \|_{\infty} < \epsilon\)} \textbf{exit} 
		      and return $u^\star = \ubarhalf$
	\While {\(
		      \psi(\ubarhalf; c, y)
		 \,{>}\,
		      \psi(\ubarnu; c, y) 
		 -
		      \nabla \psi(\ubarnu; c, y)^{\top}\!r^{\nu}
		\)
		\(
		 {}+{}
		 \tfrac{L}{2}\|r^{\nu}\|^2
		\)} (Update Lipschitz constant estimate)
	  \State Empty the L-BFGS buffers
	  \State \(L {}\gets{} 2L\), 
		 \(\sigma {}\gets{} \nicefrac{\sigma}{2}\), 
		 \(\gamma {}\gets{} \nicefrac{\gamma}{2}\)
	  \State \(\ubarhalf {}\gets{} T_{\gamma}(\ubarnu; c, y)\)
	\EndWhile	
 	\State \label{alg:panoc::d} \(d^\nu\leftarrow-H_\nu r^\nu\) using L-BFGS%
	\State \label{alg:panoc::update} 
		\(
			\ubarnuplus
		    {}\gets{}
			\ubarnu 
		    {}-{} 
			(1-\tau\!_\nu) r^\nu + \tau\!_\nu d^\nu
		\),
		where \(\tau\!_\nu\) is the largest number in \(\{\nicefrac{1}{2^i}{}:{}i\in\N\}\) such that
		\[
			\varphi_\gamma(\ubarnuplus; c, y)
		{}\leq{}
			\varphi_\gamma(\ubarnu; c, y){}-{}\sigma\|\gamma^{-1}r^\nu\|^2
		\]
\EndFor{}%
 \end{algorithmic}

\end{algorithm}

Note that Problem \eqref{eq:original_problem} can be written equivalently as
\begin{equation}
\minimize_{u\in U}\ f(u)+h(F(u)),
\end{equation}
where $h(z_1,z_2)=\delta_C(z_1)+\delta_{\{0\}}(z_2)$,
for $z_1\in\R^{n_1}$ and $z_2 \in \R^{n_2}$,
$F(u)=[F_1(u)~F_2(u)]^{\top}$, and the associated Lagrangian is
\begin{equation}
 L(u,y)=f(u)+\langle y,F(u)\rangle-h^*(y),
\end{equation}
and the first order necessary conditions for optimality can be written as
\begin{equation}
 0\in\partial_x L(u,y),\qquad 0\in\partial_y [-L(u,y)].
\end{equation}
Provided that $F_2$ is differentiable, the optimality conditions read
\begin{equation}
0\in\nabla f(u)+N_U(u)+JF(u)^\top y,\quad y\in\partial h(F(u)),
\end{equation}
or
\begin{subequations}
\begin{align}
0&\in\nabla f(u)+N_U(u)+JF_1(u)^\top y+JF_2(u)^\top y_2,\\
 y&\in\partial N_C(F(u)),\quad y_2\in\R^{n_2}.
\end{align}
\end{subequations}
It can be shown that if a pair $(\bar{u}, \bar{y})$ satisfies the termination conditions
of Algorithms \ref{alg:alm_pm} and \ref{alg:panoc}, it is an $(\epsilon, \delta)$-approximate
KKT point of Problem \eqref{eq:original_problem} in the sense that there exist
$v\in \partial_u L(\bar{u}, \bar{y})$ and $w \in \partial_y [-L(\bar u,\bar y)]$
with $\|v\| \leq \epsilon$ and $\|w\| \leq \delta$.
In particular, note that $\delta$ controls the infeasibility of the approximate solution.
The convergence properties of formulations without $F_1$ where $F_2$ is not differentiable
are studied in \citep{Hermans+2019}.

\section{Embedded Code Generation}

\subsection{Code Generation with \OpEn}
\OpEn is an open-source software for embedded nonconvex optimization which involves several independent
modules. The core numerical solver, which implements Algorithm \ref{alg:alm_pm}, is implemented in
Rust \citep{Matsakis:2014:RL:2692956.2663188}. Rust is a modern systems programming language, which
provides guarantees for memory safety, which is an important feature in embedded applications.

\OpEn provides a Python library (\opengen) and a MATLAB
toolbox for fast prototyping of embedded optimization applications; these can be used to generate
embedded Rust code for user-defined problems. \OpEn can cross-compile the generated optimizer
for several target systems (e.g., ARM Cortex-A processors). The generated
optimizer can either be used in a Rust project, or consumed through an auto-generated
\C or \Cpp interface. This allows the incorporation of the generated code in robot operating system (ROS) projects.
Furthermore, \OpEn can produce an additional interface for the auto-generated optimizer: a
very-low-latency TCP socket server, written in Rust, that uses a simple JSON data format.
The TCP server can be invoked from Python, MATLAB and virtually any programming language, while it
facilitates the deployment of edge intelligence applications \citep{Varghese+2016} and the solution of
problems by a distributed network of agents.

\begin{figure}
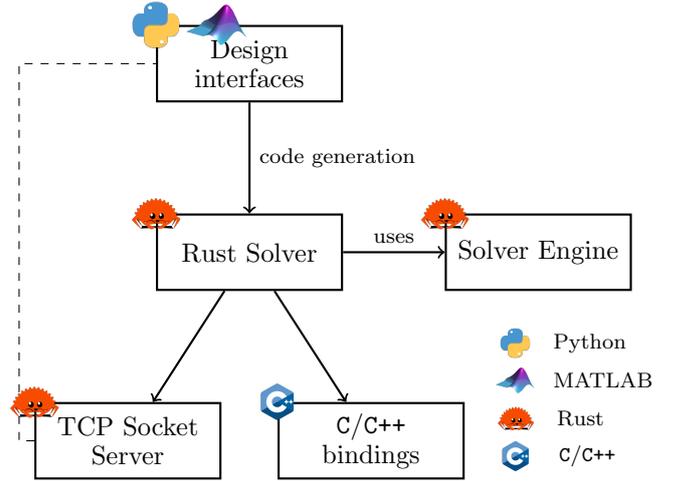

\tikzexternaldisable
{{%
	  \pgfkeys{/pgf/images/include external/.code={\includegraphics[width=0.75\columnwidth]{#width=0.75\columnwidth}}}%
	  \tikzsetnextfilename{OpEn}%
	  \input{./Tikz/OpEn.tex}%
    }}
\caption{Structure of the prototyping and embedded code generation framework of \OpEn.}
\end{figure}

In order to use \OpEn in Python, the user only needs to install the Rust compiler,
\texttt{clang} and \opengen using \texttt{pip}:
\begin{lstlisting}[style=myPython]
pip install opengen
\end{lstlisting}

To use \opengen in Python, one first needs to import it alongside CasADi (which is
co-installed with \opengen)
\begin{lstlisting}[style=myPython]
import opengen as og
import opengen.constraints as ogc
import casadi.casadi as cs
\end{lstlisting}

Let us give a simple example of using \OpEn to generate a parametric optimizer.
Consider the following problem of constrained minimization of the Rosenbrock
function:
\begin{subequations}
\begin{align}
 \mathbb{P}(p){}:{}\minimize_{u\in \R^{5}}&  \sum_{i=1}^{4}p_2(u_{i+1}-u_i^2)^2 + (p_1-u_i)^2
 \\
 \subjto& \|u\| \leq 0.73
 \\
 &p_3\sin(u_1) = \cos(u_2 + u_3)
 \\
 &u_3 + u_4 \leq 0.2
\end{align}
\end{subequations}
with $p\in\R^3$.
This is a nonconvex parametric optimization problem which can be written in
multiple ways in the form of Problem \eqref{eq:original_problem}. One possible
choice is $U=\{u\in\R^5 {}:{} \|u\| \leq 0.73\}$ and $F_2:\R^5\times\R^3\to\R^2$
is given by
\begin{equation}
    F_2(u, p)
 {}={}
    \begin{bmatrix}
      p_3\sin(u_1) - \cos(u_2 + u_3)
      \\
      \plus{u_3 + u_4 - 0.2}
    \end{bmatrix}.
\end{equation}
In this case, both constraints will be handled with the penalty method.

\begin{lstlisting}[style=myPython]
u = cs.SX.sym("u", 5) # decision variable
p = cs.SX.sym("p", 3) # parameter

# Define cost function
f = sum([p[1] * (u[i+1] - u[i]**2)**2
         + (p[0] - u[i])**2 for i in range(4)])

# Define constraints
c1 = p[2]*cs.sin(u[0]) - cs.cos(u[1] + u[2])
c2 = cs.fmax(u[2] + u[3] - 0.2, 0)
c = cs.vertcat(c1, c2)
bounds = ogc.Ball2(radius=0.73)

# Formulate problem
problem = og.builder.Problem(u, p, f) \
    .with_penalty_constraints(c)      \
    .with_constraints(bounds)

# Configure the build and the solver
build_cfg = og.config.BuildConfiguration() \
    .with_build_directory("build")         \
    .with_tcp_interface_config()           \
    .with_build_c_bindings()
meta = og.config.OptimizerMeta() \
    .with_optimizer_name("my_optimizer")
solver_cfg = og.config.SolverConfiguration() \
    .with_tolerance(1e-5)                    \
    .with_delta_tolerance(1e-4)              \
    .with_initial_tolerance(1e-4)            \
    .with_initial_penalty(1e3)               \
    .with_penalty_weight_update_factor(5)

# Generate code
builder = og.builder.OpEnOptimizerBuilder(
    problem, meta, build_cfg, solver_cfg)
builder.build()
\end{lstlisting}

A second option is to handle both constraints with the augmented
Lagrangian method by defining
\begin{equation}
 F_1(u, p) = \begin{bmatrix}
      p_3\sin(u_1) - \cos(u_2 + u_3)
      \\
      u_3 + u_4 - 0.2
    \end{bmatrix}
\end{equation}
and $C = \{0\} \times (-\infty, 0]$. In that case, set $Y$ is computed internally and
is of the form $Y = [-M, M] \times [0, M]$, where $M=10^{12}$ (however, the user can
override this and provide her own set compact $Y$).

\begin{lstlisting}[style=myPython]
# Set C = {0} x (-inf, 0]
n1 = 2
Ca = ogc.Zero()  # set {0}
Cb = ogc.Rectangle(None, [0])  # set (-inf, 0]
C = ogc.CartesianProduct(n1, [0, 1], [Ca, Cb])

# Custom set Y (Optional)
M = 1e10
Y = ogc.Rectangle([-M, 0], [M, M])

# Problem specification
problem = og.builder.Problem(u, p, f)\
    .with_aug_lagrangian_constraints(c, C, Y)\
    .with_constraints(bounds)
\end{lstlisting}

The solver will be generated and all relevant files will be stored in
\texttt{build/my\_optimizer}. The solver can be interfaced through its \C/\Cpp
interface (enabled with \texttt{.with\_build\_c\_bindings}) using the
auto-generated bindings \texttt{build/my\_optimizer/my\_optimizer\_bindings.h}.
It can also be invoked over a generated TCP socket interface which
is enabled using \texttt{with\_tcp\_interface\_config} (the default IP and
port are \texttt{127.0.0.1} and \texttt{8333}, but this can be configured).

\begin{lstlisting}[style=myPython]
# Create a TCP manager
optimizer_path = 'build/my_optimizer'
mng = og.tcp.OptimizerTcpManager(optimizer_path)

# Start the server
mng.start()

# Use the solver (returns the inner/outer iter-
# ations, solution time, infeasibility, solution
# status, Lagrange multipliers and solution)
solver_status = mng.call([1.0, 50.0, 1.5])

# Kill the server
mng.kill()
\end{lstlisting}

For $p=(1,50,1.5)$, the formulation with the penalty method led to 7 outer iterations,
647 inner iterations in total and executed in $\unit[3.5]{ms}$. The formulation based on
the augmented Lagrangian method ran in $\unit[1.4]{ms}$ after 5 outer and 175 total inner
iterations. We will discuss the convergence speed of the proposed method in Section
\ref{sec:simulations} in greater detail. IPOPT solves the same problem in $\unit[8.2]{ms}$
and \texttt{scipy}'s SQP solver, SLSQP, which is based on \citep{KraftSLSQP88}, in $\unit[15.3]{ms}$.

\subsection{Implementation Details}
While the solver itself is written in Rust, the generation of code for the cost function, its gradient,
and functions $F_1$ and $F_2$ is delegated to CasADi \citep{Andersson2019} which generates \C code.
This optimized \C code is compiled into the Rust binary using the Rust tools \emph{bindgen} and
\emph{cc}, which gives a Rust interface to the generated code.
The generation of both the Rust solver and the TCP socket server is facilitated by the
\emph{Jinja2 template language}.
\C/\Cpp, as well as any language that supports the \C/\Cpp application binary interface (ABI),
can consume the generated optimizers. To that end, the Rust tool \emph{cbindgen} is used to
generate appropriate bindings. This allows the solver to be readily integrated with the robot
operating system \citep{Koubaa:2017:ROS:3131501}.
These design choices allow the entire solver and library to run in Rust thus retaining its strong
memory safety guarantees.

Moreover, the solver is thoroughly and continuously being tested in \emph{Travis CI}
(108 Rust tests, 35 Python tests) with 96.3\% code coverage, alongside extensive
documentation of both \OpEn's Rust code and \opengen's Python code, and a
large collection of examples available on the website at \texttt{http://doc.optimization-engine.xyz}.

\section{Simulations}\label{sec:simulations}

\subsection{Obstacle Avoidance}
The dynamics of an autonomous ground vehicle can be described by the following nonlinear
bicycle model with four state and two input variables
\begin{subequations}
 \begin{align}
  \dot{p}_x {}={}& v \cos \psi,
  \\
  \dot{p}_y {}={}& v \sin \psi,
  \\
  \dot{\psi} {}={}& \nicefrac{v}{L}\tan \delta,
  \\
  \dot{v} {}={}& \alpha (a - v),
 \end{align}
\end{subequations}
where $(p_x, p_y)$ are the position coordinates of the vehicle with respect to an
Earth-fixed frame of reference, $\psi$ is the vehicle's orientation and $v$ is the
longitudinal velocity, which follows a first-order dynamics with parameter $\alpha = 0.25$.
The vehicle is controlled by its longitudinal acceleration, $a$, and the steering angle $\delta$.
We define the state vector $x=(p_x, p_t, \psi, v)$ and the input vector $u=(a, \delta)$.
The length of the vehicle is $L = \unit[0.5]{m}$.
The above dynamical system is discretized with the Euler method with sampling time
$T_s=\unit[50]{ms}$.
The vehicle must navigate from an initial position, pose and velocity, $x$, to a target position
and pose,
\(
      x^{\mathrm{ref}}
      =
      (p_x^{\mathrm{ref}},
       p_y^{\mathrm{ref}},
       \psi^{\mathrm{ref}},
       0)
\),
while avoiding a cylindrical obstacle of radius $r = 0.65$ positioned at $(-3, 0.2)$.

\begin{figure}
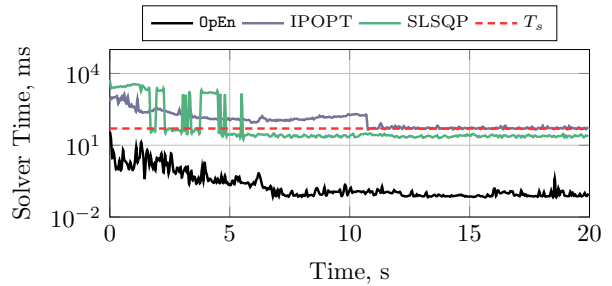

 \centering
 {{%
	  \pgfkeys{/pgf/images/include external/.code={\includegraphics[width=0.9\columnwidth]{#width=0.9\columnwidth}}}%
	  \tikzsetnextfilename{obstacleAndPath_times}%
	  \input{./Tikz/obstacleAndPath_times.tex}%
    }}
 \caption{Computation time for the NMPC controlling the autonomous ground vehicle.
 (Solid lines) Execution times for \OpEn, IPOPT and SLSQP, (Dashed line) sampling time.}
\label{fig:obstacle_times_comparison}
\end{figure}

We formulate an MPC problem with prediction horizon $N=100$ and stage cost function
$\ell(x, u)=18(p_x^2+p_y^2)+2\psi^2+5v^2$ and the terminal cost function $\ell_N(x)=
1500(p_x^2+p_y^2)+500\psi^2+10v^2$. In order to avoid aggressive steering or jerks we
introduce the additional stage cost function
\(
      \ell_\Delta(u_t, u_{t-1})
  {}={}
      100(a_{t}-a_{t-1})^2 + 30(\delta_t-\delta_{t-1})^2
\).
Furthermore, we impose the actuation constraints $-1 \leq a \leq \unitfrac[2]{m}{s^2}$
and $|\delta| \leq \unit[0.25]{rad}\approx 14.32^\circ$.
Then, the state sequence is eliminated as in Equation \eqref{eq:single_shooting}.
The obstacle avoidance constraints can be formulated as follows
\begin{equation}\label{eq:obstacle:avoidance_constraints}
 \plus{r^2 - p_x^2 - p_y^2} = 0,
\end{equation}
which can be accommodated using $F_2$ and the penalty method. Alternatively, in 
this case we can treat the obstacle constraints as inequality constraints of 
the form $r^2 - p_x^2 - p_y^2 \leq 0$ and describe them using $F_1$. The two 
approaches lead to comparable results as shown in Figure~\ref{fig:obstacle_times_comparison}.

\tikzexternalenable
\begin{figure}[ht]
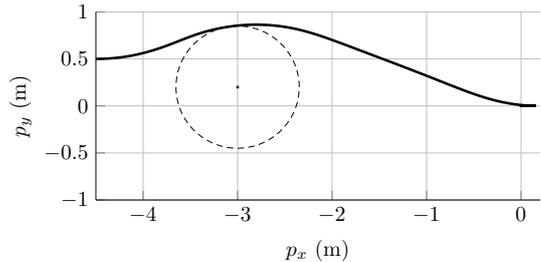

 \centering
 {{%
	  \pgfkeys{/pgf/images/include external/.code={\includegraphics[width=0.8\columnwidth]{#width=0.8\columnwidth}}}%
	  \tikzsetnextfilename{obstacleAndPath}%
	  \input{./Tikz/obstacleAndPath.tex}%
    }}
 \caption{(Solid line) Path of NMPC-controlled ground vehicle. (Dashed line) Circular
          obstacle at $(x_c, y_c) = (-3, 0.2)$ of radius $r=\unit[0.65]{m}$.}
 \label{fig:obstacle:2dpath}           
\end{figure}

\begin{figure}
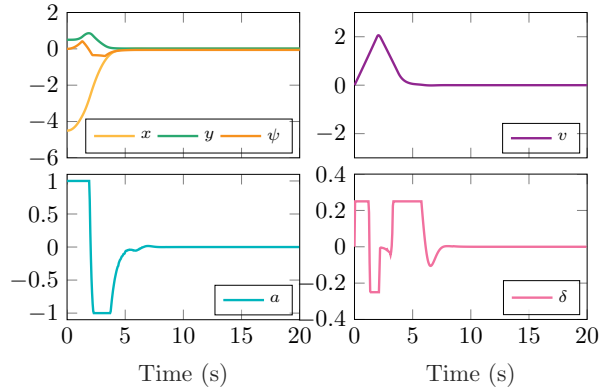

 \centering
 {{%
	  \pgfkeys{/pgf/images/include external/.code={\includegraphics[width=0.9\columnwidth]{#width=0.9\columnwidth}}}%
	  \tikzsetnextfilename{obstacleAndPath_traj}%
	  \input{./Tikz/obstacleAndPath_traj.tex}%
    }}
 \caption{Closed-loop trajectories for the NMPC-controlled vehicle.}
 \label{fig:obstacle:path_time_profile}
\end{figure}

The NMPC-controlled autonomous vehicle manages to avoid the obstacle, arrive
at the target position and assume the desired orientation as shown in Figures
\ref{fig:obstacle:2dpath} and \ref{fig:obstacle:path_time_profile}.

\OpEn outperforms significantly the interior point solver IPOPT and the SQP solver
of Python's package \texttt{scipy}, SLSQP, as illustrated in Figure
\ref{fig:obstacle_times_comparison}. All software make use of the same auto-generated
code to evaluate the cost function, $f$, its gradient, $\nabla f$, and function $F_2$
and the tolerance is set to $\delta=10^{-3}$ and $\epsilon=10^{-4}$.
The selected solver parameters are
$\epsilon_0=10^{-4}$,
$\rho=5$,
$c_0=500$ and L-BFGS memory
$\mu=20$.

\subsection{Constrained Nonlinear Estimation}
Consider Lorenz's chaotic system
\begin{equation}
 \dot{x} =
\begin{bmatrix}
a_1 (x_{2} - x_{1})
\\
x_{1} (a_2 - x_{3}) - x_{2}
\\
x_{1} x_{2} - a_3 x_{3}
\end{bmatrix} + w(t),
\end{equation}
where $w$ is a disturbance signal, and parameters $a_1=10$, $a_2=14$, $a_3=\nicefrac{8}{3}$, which we discretize by means
of the fourth-order Runge-Kutta method with integration step $h = 0.1$ leading to a discrete-time system
$x_{t+1} = \Phi(x_t) + w_t$.
It is known that $-1 \leq w \leq 1$.
The system output is given by
\begin{equation}
 y =
\smashunderbracket{\begin{bmatrix}
2x_{1}
\\
x_{2} + x_{3}
\end{bmatrix}}{G(x)} + v,
\vspace{1.15em}
\end{equation}
where $v$ is a noise signal for which it is known that $-1.5 \leq v \leq 1.5$.
We formulate the following constrained nonconvex estimation problem:
\begin{subequations}\label{eq:nmhe:problem_statement}
\begin{align}
\mathbb{P}(\bm{y}){}:{}\minimize_{\hat{\bm{x}}, \hat{\bm{w}}, \hat{\bm{v}}} &
    \sum_{t=0}^{N-1}\|\hat{w}_t\|^2_Q + \|\hat{v}_t\|^2_R
    \label{eq:nmhe:cost}
    \\
    \subjto\,& w_{\mathrm{min}} \leq \hat{w}_{t} \leq w_{\mathrm{max}}, t\in\N_{N}
    \label{eq:nmhe:wconstraints}
    \\
    &v_{\mathrm{min}} \leq \hat{v}_{t} \leq v_{\mathrm{max}}, t\in\N_{N}
    \label{eq:nmhe:vconstraints}
    \\
    &\hat{x}_{t+1} = \Phi(\hat{x}_t) + \hat{w}_t, t\in\N_{N-1}
    \label{eq:nmhe:dynamics}
    \\
    & y_t = G(\hat{x}_t) + \hat{v}_t, t\in\N_{N}
    \label{eq:nmhe:ouput}
\end{align}
\end{subequations}
where $Q$ and $R$ are symmetric positive definite matrices.
This is a problem with decision variables $\hat{\bm{x}} = (\hat{x}_t)_{t=0}^{N}$,
\(\hat{\bm{w}}=(\hat{w}_t)_{t=0}^{N}\)
and
\(\hat{\bm{v}}=(\hat{v}_t)_{t=0}^{N}\),
that is, $u=(\hat{\bm{x}}, \hat{\bm{w}}, \hat{\bm{v}})$,
and parameter $\bm{y}=(y_t)_{t=0}^{N}$.

Constraints \eqref{eq:nmhe:wconstraints} and \eqref{eq:nmhe:vconstraints} define a rectangle,
on which one can easily compute projections. The equality constraints \eqref{eq:nmhe:dynamics}
and \eqref{eq:nmhe:ouput} are treated with the augmented Lagrangian method by defining
\begin{equation}
 F_1(u, \bm{y}) =
 \begin{bmatrix}
  (\hat{x}_{t+1} - \Phi(\hat{x}_t) - \hat{w}_t)_{t=0}^{N-1}
  \\
  (y_t - G(\hat{x}_t) - \hat{v}_t)_{t=0}^{N}
 \end{bmatrix},
\end{equation}
and $C=\{0\}$.

\begin{figure}
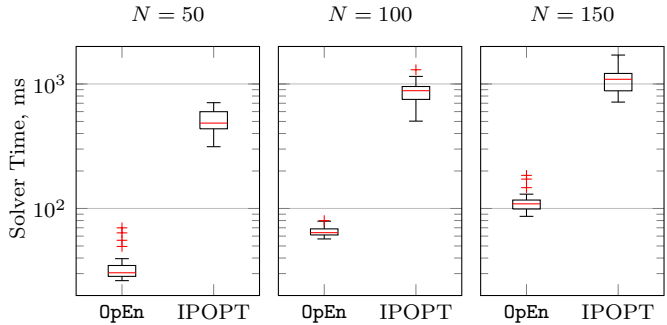

 \centering
 {{%
	  \pgfkeys{/pgf/images/include external/.code={\includegraphics[]{#}}}%
	  \tikzsetnextfilename{nmhe_times_50}%
%
%
\begin{tikzpicture}

\begin{axis}[%
width=0.95in,
height=1.3in,
at={(0.772in,0.484in)},
scale only axis,
xmin=0.5,
xmax=2.5,
xtick={1,2},
xticklabels={{\footnotesize \OpEn},{\footnotesize IPOPT}},
ymin=20,
ymax=2000,
axis background/.style={fill=white},
ylabel={{\footnotesize Solver Time, ms}},
y label style={at={(0.2,0.55)}},
ymode=log,
title={\footnotesize $N=50$},
ticklabel style = {font=\footnotesize},
ymajorgrids
]
\addplot [color=black,  forget plot]
  table[row sep=crcr]{%
1	34.851627\\
1	39.544335\\
};
\addplot [color=black,  forget plot]
  table[row sep=crcr]{%
2	598.888874053955\\
2	707.505702972412\\
};
\addplot [color=black,  forget plot]
  table[row sep=crcr]{%
1	26.338626\\
1	28.501168\\
};
\addplot [color=black,  forget plot]
  table[row sep=crcr]{%
2	313.437461853027\\
2	436.547517776489\\
};
\addplot [color=black, forget plot]
  table[row sep=crcr]{%
0.925	39.544335\\
1.075	39.544335\\
};
\addplot [color=black, forget plot]
  table[row sep=crcr]{%
1.925	707.505702972412\\
2.075	707.505702972412\\
};
\addplot [color=black, forget plot]
  table[row sep=crcr]{%
0.925	26.338626\\
1.075	26.338626\\
};
\addplot [color=black, forget plot]
  table[row sep=crcr]{%
1.925	313.437461853027\\
2.075	313.437461853027\\
};
\addplot [color=black, forget plot]
  table[row sep=crcr]{%
0.85	28.501168\\
0.85	34.851627\\
1.15	34.851627\\
1.15	28.501168\\
0.85	28.501168\\
};
\addplot [color=black, forget plot]
  table[row sep=crcr]{%
1.85	436.547517776489\\
1.85	598.888874053955\\
2.15	598.888874053955\\
2.15	436.547517776489\\
1.85	436.547517776489\\
};
\addplot [color=red, forget plot]
  table[row sep=crcr]{%
0.85	30.4637095\\
1.15	30.4637095\\
};
\addplot [color=red, forget plot]
  table[row sep=crcr]{%
1.85	484.921097755432\\
2.15	484.921097755432\\
};
\addplot [color=black, draw=none, mark=+, mark options={solid, red}, forget plot]
  table[row sep=crcr]{%
1	49.548377\\
1	55.56271\\
1	63.72142\\
1	69.950045\\
};
\end{axis}
\end{tikzpicture}%
%
    }}%
 {{%
	  \pgfkeys{/pgf/images/include external/.code={\includegraphics[]{#}}}%
	  \tikzsetnextfilename{nmhe_times_100}%
%
%
\begin{tikzpicture}

\begin{axis}[%
width=0.95in,
height=1.3in,
at={(0.772in,0.484in)},
scale only axis,
xmin=0.5,
xmax=2.5,
xtick={1,2},
xticklabels={{\footnotesize \OpEn},{\footnotesize IPOPT}},
yticklabels={,,},
ymin=20,
ymax=2000,
axis background/.style={fill=white},
ymode=log,
title={\footnotesize $N=100$},
ymajorgrids,
]
\addplot [color=black,  forget plot]
  table[row sep=crcr]{%
1	68.595252\\
1	79.086087\\
};
\addplot [color=black,  forget plot]
  table[row sep=crcr]{%
2	954.029560089111\\
2	1149.82318878174\\
};
\addplot [color=black,  forget plot]
  table[row sep=crcr]{%
1	56.954211\\
1	61.297961\\
};
\addplot [color=black,  forget plot]
  table[row sep=crcr]{%
2	502.937793731689\\
2	751.757144927979\\
};
\addplot [color=black, forget plot]
  table[row sep=crcr]{%
0.925	79.086087\\
1.075	79.086087\\
};
\addplot [color=black, forget plot]
  table[row sep=crcr]{%
1.925	1149.82318878174\\
2.075	1149.82318878174\\
};
\addplot [color=black, forget plot]
  table[row sep=crcr]{%
0.925	56.954211\\
1.075	56.954211\\
};
\addplot [color=black, forget plot]
  table[row sep=crcr]{%
1.925	502.937793731689\\
2.075	502.937793731689\\
};
\addplot [color=black, forget plot]
  table[row sep=crcr]{%
0.85	61.297961\\
0.85	68.595252\\
1.15	68.595252\\
1.15	61.297961\\
0.85	61.297961\\
};
\addplot [color=black, forget plot]
  table[row sep=crcr]{%
1.85	751.757144927979\\
1.85	954.029560089111\\
2.15	954.029560089111\\
2.15	751.757144927979\\
1.85	751.757144927979\\
};
\addplot [color=red, forget plot]
  table[row sep=crcr]{%
0.85	63.9673985\\
1.15	63.9673985\\
};
\addplot [color=red, forget plot]
  table[row sep=crcr]{%
1.85	882.893443107605\\
2.15	882.893443107605\\
};
\addplot [color=black, draw=none, mark=+, mark options={solid, red}, forget plot]
  table[row sep=crcr]{%
1	80.011753\\
};
\addplot [color=black, draw=none, mark=+, mark options={solid, red}, forget plot]
  table[row sep=crcr]{%
2	1300.4994392395\\
};
\end{axis}
\end{tikzpicture}%
%
    }}%
 {{%
	  \pgfkeys{/pgf/images/include external/.code={\includegraphics[]{#}}}%
	  \tikzsetnextfilename{nmhe_times_150}%
%
%
\begin{tikzpicture}

\begin{axis}[%
width=0.95in,
height=1.3in,
at={(0.772in,0.484in)},
scale only axis,
xmin=0.5,
xmax=2.5,
xtick={1,2},
xticklabels={{\footnotesize \OpEn},{\footnotesize IPOPT}},
yticklabels={,,},
ymin=20,
ymax=2000,
axis background/.style={fill=white},
ymode=log,
title={\footnotesize $N=150$},
ymajorgrids
]
\addplot [color=black,  forget plot]
  table[row sep=crcr]{%
1	116.895546\\
1	130.331505\\
};
\addplot [color=black,  forget plot]
  table[row sep=crcr]{%
2	1214.80703353882\\
2	1706.72345161438\\
};
\addplot [color=black,  forget plot]
  table[row sep=crcr]{%
1	86.473378\\
1	99.053962\\
};
\addplot [color=black,  forget plot]
  table[row sep=crcr]{%
2	715.24977684021\\
2	881.118535995483\\
};
\addplot [color=black, forget plot]
  table[row sep=crcr]{%
0.925	130.331505\\
1.075	130.331505\\
};
\addplot [color=black, forget plot]
  table[row sep=crcr]{%
1.925	1706.72345161438\\
2.075	1706.72345161438\\
};
\addplot [color=black, forget plot]
  table[row sep=crcr]{%
0.925	86.473378\\
1.075	86.473378\\
};
\addplot [color=black, forget plot]
  table[row sep=crcr]{%
1.925	715.24977684021\\
2.075	715.24977684021\\
};
\addplot [color=black, forget plot]
  table[row sep=crcr]{%
0.85	99.053962\\
0.85	116.895546\\
1.15	116.895546\\
1.15	99.053962\\
0.85	99.053962\\
};
\addplot [color=black, forget plot]
  table[row sep=crcr]{%
1.85	881.118535995483\\
1.85	1214.80703353882\\
2.15	1214.80703353882\\
2.15	881.118535995483\\
1.85	881.118535995483\\
};
\addplot [color=red, forget plot]
  table[row sep=crcr]{%
0.85	108.988067\\
1.15	108.988067\\
};
\addplot [color=red, forget plot]
  table[row sep=crcr]{%
1.85	1090.24012088776\\
2.15	1090.24012088776\\
};
\addplot [color=black, draw=none, mark=+, mark options={solid, red}, forget plot]
  table[row sep=crcr]{%
1	147.025007\\
1	172.214507\\
1	184.441799\\
};
\end{axis}
\end{tikzpicture}%
%
    }}%
 \caption{Runtime of \OpEn and IPOPT; box plots for 30 random trials.}
 \label{fig:nmhe:times_boxplot}
\end{figure}

\tikzexternaldisable
\begin{figure}
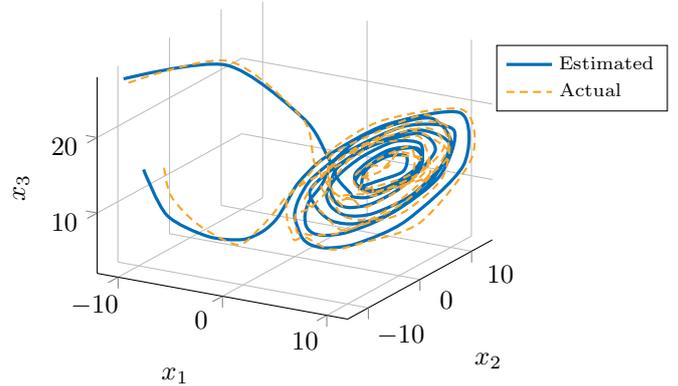

 \centering
 {{%
	  \pgfkeys{/pgf/images/include external/.code={\includegraphics[width=0.75\columnwidth]{#width=0.75\columnwidth}}}%
	  \tikzsetnextfilename{Estimation}%
	  \input{./Tikz/Estimation.tex}%
    }}
 \caption{Constrained state estimation for Lorenz's system with $N=100$}
 \label{fig:nmhe:estimated_vs_actual}
\end{figure}

\begin{figure}
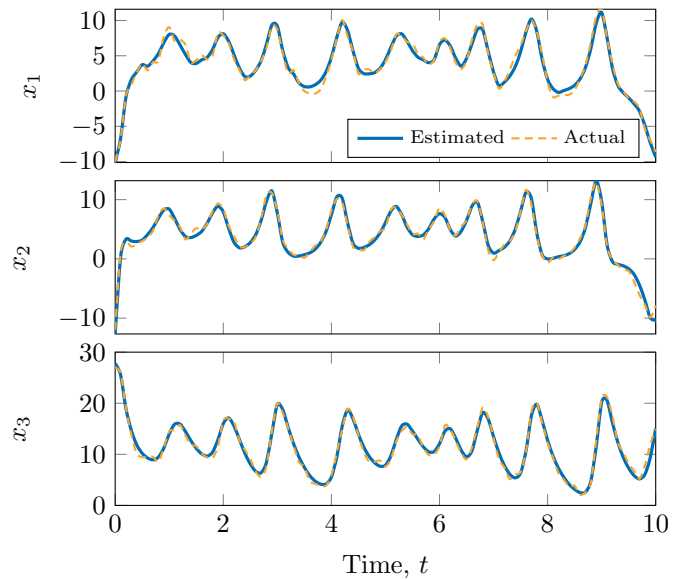

 \centering
 {{%
	  \pgfkeys{/pgf/images/include external/.code={\includegraphics[width=0.75\columnwidth]{#width=0.75\columnwidth}}}%
	  \tikzsetnextfilename{x3_estimation}%
	  \input{./Tikz/x3_estimation.tex}%
    }}
 \caption{States and their estimates for Lorenz's system 
          using the constrained nonlinear estimation formulation 
          given in \eqref{eq:nmhe:problem_statement} with $N=100$}
          \label{fig:nmhe:estimated_vs_actual_time}
\end{figure}

In Figure \ref{fig:nmhe:times_boxplot} we compare the runtime of \OpEn with IPOPT
for $N=50$, $N=100$ and $N=150$.
The SQP solver of \texttt{scipy} was significantly slower compared to both \OpEn and
IPOPT and was therefore omitted.
An estimated trajectory for the case $N=100$ is shown in Figure \ref{fig:nmhe:estimated_vs_actual}.
The same data are presented in Figure \ref{fig:nmhe:estimated_vs_actual_time} where the
system states and their estimates are plotted against time.
The tuning parameters for \OpEn were chosen to be $c_0=200$, $\rho=1.8$, $\epsilon_0=0.1$,
$\mu=15$ and tolerances $\delta=10^{-5}$ and $\epsilon=10^{-4}$.
\OpEn performed no more that 7 outer iterations (in the majority of cases, fewer than 5)
keeping the penalty parameter, $c_\nu$, below $39672$.

\section{Conclusions and Future work}
In this paper we presented a code generation software that can be easily used 
in Python and MATLAB and will empower engineers to build embedded optimization 
modules that will underpin complex control, estimation and signal processing 
developments. The proposed algorithm can handle nonconvex smooth cost functions
and nonlinear constraints by means of the augmented Lagrangian and penalty methods,
the latter being particularly suitable for obstacle/collision avoidance problems.
Future work will focus on the study of the convergence properties of problems 
that involve constraints of the form \eqref{eq:constraints:alm} and \eqref{eq:constraints:pm}
under weak assumptions.

\bibliography{ifacconf}

\end{document}